\documentclass[12pt,leqno,draft]{article} 
\usepackage{amsmath,amssymb,amsthm,amsfonts}
\usepackage{enumerate,color}
\makeatletter
\def\@cite#1#2{[{{\bfseries #1}\if@tempswa , #2\fi}]}
\renewcommand{\section}{%
\@startsection{section}{1}{\z@}
{0.5truecm plus -1ex minus -.2ex}%
{1.0ex plus .2ex}{\bfseries\large}}
\def\@seccntformat#1{\csname the#1\endcsname.\ }
\makeatother

\numberwithin{equation}{section} 
\newtheorem{thm}{Theorem}[section]

\newtheorem{prop}[thm]{Proposition}
\theoremstyle{definition}

\newtheorem{remark}{Remark}[section]

\begin{document}
\begin{center}
    \Large{{\bf 
             Remarks on two connected papers about Keller--Segel systems with nonlinear production
          }}
\end{center}
\vspace{5pt}
\begin{center}
    Yuya Tanaka\footnote{Corresponding author.}\\
               \vspace{8pt}
    \footnote[0]{{E-mail addresses}:
    {\tt yuya.tns.6308@gmail.com}, 
    {\tt giuseppe.viglialoro@unica.it},\\ \qquad
    {\tt yokota@rs.tus.ac.jp}}
    Department of Mathematics,
    Tokyo University of Science\\
    1-3, Kagurazaka, Shinjuku-ku, 
    Tokyo 162-8601, Japan\\
               \vspace{15pt}
    Giuseppe Viglialoro\\
               \vspace{8pt}
    Dipartimento di Matematica e Informatica, 
    Università di Cagliari,\\
    Via Ospedale 72, 09124,
    Cagliari, Italy\\
               \vspace{15pt}
    Tomomi Yokota\\
               \vspace{8pt}
    Department of Mathematics,
    Tokyo University of Science\\
    1-3, Kagurazaka, Shinjuku-ku, 
    Tokyo 162-8601, Japan\\
\end{center}
\begin{center}    
    \small \today
\end{center}
\vspace{2pt}
\newenvironment{summary}
{\vspace{.5\baselineskip}\begin{list}{}{%
     \setlength{\baselineskip}{0.85\baselineskip}
     \setlength{\topsep}{0pt}
     \setlength{\leftmargin}{12mm}
     \setlength{\rightmargin}{12mm}
     \setlength{\listparindent}{0mm}
     \setlength{\itemindent}{\listparindent}
     \setlength{\parsep}{0pt}
     \item\relax}}{\end{list}\vspace{.5\baselineskip}}
\begin{summary}
{\footnotesize {\bf Abstract.}
These notes aim to provide a deeper insight on the specifics of two articles dealing with chemotaxis models with nonlinear production. More precisely, we are referring to the papers ``Boundedness of solutions to a quasilinear parabolic--parabolic chemotaxis model with nonlinear signal production'' by X. Tao, S. Zhou and M. Ding [\textit{J. Math.\ Anal.\ Appl.} {\bf 474}:1 (2019) 733--747] and  ``Boundedness for a fully parabolic {K}eller--{S}egel model	with sublinear segregation and superlinear aggregation'' by S. Frassu and G. Viglialoro [\textit{Acta Appl.\ Math.}
	{\bf 171}:1 (2021), 19]. These works, independently published in these last years, present results leaving open room for further improvement. Indeed, in the first a gap in the proof of the main claim appears, whereas the cornerstone assumption in the second is not sharp. In these pages we give a more complete picture to the relative underlying comprehension.
	
%
}
\end{summary}
\newpage
%
%
\section{Motivations and main result}
%
%
In this short document we focus on \cite[Theorem 1.1]{TZD19} and \cite[Theorem 2.1]{FV21} where chemotaxis models for two coupled parabolic equations are so formulated:
  \begin{align}\label{P}
    \begin{cases}
    u_t=\nabla \cdot (D(u) \nabla u) - \nabla \cdot (S(u) \nabla v), &x \in \Omega,\ t>0,\\
    v_t=\Delta v - v + g(u), &x \in \Omega,\ t>0,\\
    \frac{\partial u}{\partial \nu} = \frac{\partial v}{\partial \nu} =0, &x \in \partial\Omega,\  t>0,\\
    u(x,0)=u_0(x), \quad v(x,0)=v_0(x), &x \in \Omega,\ t>0.
    \end{cases}
  \end{align}
Herein, $\Omega \subset \mathbb{R}^n \ (n\ge2)$ is a bounded domain with smooth boundary, and $\frac{\partial}{\partial \nu}$ denotes the differentiation with respect to the outward normal of $\partial \Omega$.  Additionally, the initial data $(u_0,v_0)$ is assumed to satisfy
  \begin{align}\label{initial}
    \begin{cases}
    u_0 \in C^0(\overline{\Omega}) \quad\mbox{is nonnegative with}\ u_0 \not\equiv 0,\\
    v_0 \in C^1(\overline{\Omega}) \quad\mbox{is nonnegative},
    \end{cases}
  \end{align}
whereas, for all $u\geq 0$ and appropriate  real numbers $d_0, d_1, s_1,\alpha, \alpha_1, \beta, g_1, \gamma$, the diffusion and sensitivity laws  $D,S \in C^2([0,\infty))$ and the production growth 
$g \in C^1([0,\infty))$  are such that
  \begin{align}\label{DS}
    d_0(1+u)^{-\alpha} \le D(u) \le d_1(1+u)^{-\alpha_1},\quad
    0 \le S(u) \le s_1u(1+u)^{\beta-1},
  \end{align}
and 
  \begin{align}\label{g}
    0 \le g(u) \le g_1u^\gamma.
  \end{align}
The aforementioned results in \cite{TZD19} and \cite{FV21} are collected as follows.
\begin{thm}
Let $n \ge 2$ and $(u_0,v_0)$ satisfy \eqref{initial}.
Suppose that $D, S$ and $g$ fulfill \eqref{DS} and \eqref{g}. Then problem \eqref{P} admits a unique nonnegative classical solution $(u,v)$ which is globally bounded provided that:
\begin{enumerate}[I)]
\item {\cite[Theorem 1.1]{TZD19}}\label{thm} $0<\gamma \le 1$ and 
\begin{align}\label{TZDcondi}
	\alpha + \beta + \gamma <1 + \frac{2}{n};
\end{align}
\item {\cite[Theorem 2.1]{FV21}}\label{thmFraVig} $\alpha=\alpha_1=0$,  $0<\gamma <\frac{2}{n}$, $\beta \ge \frac{2}{n}$ and 
\begin{align}\label{FraVigcondi}
	\beta + \frac{\gamma}{2} <1 + \frac{1}{n}.
\end{align}

\end{enumerate}
\end{thm}
These two theorems have been  proved, in an independent way the one from the other, recently. Moreover,  when investigating a variant of Keller--Segel systems like those in \eqref{P}, the authors of this report realized that: 
\begin{itemize}
	\item  for $0<\gamma< \frac{1}{n}$, the proof leading to condition \eqref{TZDcondi} has a mathematical inconsistency; in this same range, even for the linear diffusion case $\alpha=\alpha_1=0$, the condition cannot hold true and has to be replaced by \eqref{FraVigcondi}; 
	\item for $\frac{1}{n}\leq \gamma<\frac{2}{n}$ and $\alpha=\alpha_1=0$, assumption  \eqref{FraVigcondi} is less accurate than \eqref{TZDcondi}.
\end{itemize}
Since this gap leaves the general theory about models \eqref{P} somehow incomplete and fragmented,  we understand that it is of primary importance giving a revised and unified conclusion. Precisely, the role behind the forthcoming theorem is twofold:  correcting  \cite[Theorem 1.1]{TZD19} and improving \cite[Theorem 2.1]{FV21}. 
\begin{thm}\label{Mainthm}
	Let $n \ge 2$ and $(u_0,v_0)$ satisfy \eqref{initial}.
	Suppose that $D, S$ and $g$ fulfill \eqref{DS} and \eqref{g}.
	If $0<\gamma \le 1$ and 
  \begin{align}\label{newcondi}
	\begin{cases}
		\alpha + \beta + \gamma<1 + \frac{2}{n} \quad&\mbox{if }\ \gamma \in \left[\frac{1}{n},1\right],\\
		\alpha + \beta <1 + \frac{1}{n} \quad&\mbox{if }\ \gamma \in \left(0,\frac{1}{n}\right),
	\end{cases}
\end{align}
	then problem \eqref{P} admits a unique nonnegative classical solution $(u,v)$ which is globally bounded.
\end{thm}
\section{Identification of the gap}
Once combined with well-known extensibility criteria, global boundedness for local classical solutions to problem \eqref{P}, defined in $\Omega \times (0,T_{{\rm max}})$, is achieved by controlling $ \|u(\cdot,t)\|_{L^p(\Omega)}$ and $
\|\nabla v(\cdot,t)\|_{L^q(\Omega)} $ on $(0,T_{\rm max})$, and for $p,q$ large enough. In particular, if we refer to \cite{TZD19}, such boundedness  relies on the ensuing
\begin{prop}[{\cite[Proposition 3.1]{TZD19}}]\label{prop}
Let $n \ge 2$ and $(u_0,v_0)$ satisfy \eqref{initial}. 
Suppose that $D, S$ and $g$ fulfill \eqref{DS} and \eqref{g}.
If $0<\gamma \le1$, $\alpha$ and $\beta$ are constrained by assumption \eqref{TZDcondi}, then for all $p \in [1,\infty)$ and each $q \in [1,\infty)$,
there exists $C=C(p,q,\alpha,\alpha_1,\beta,\gamma)>0$ such that 
  \[
    \|u(\cdot,t)\|_{L^p(\Omega)} \le C \quad\mbox{and}\quad
    \|\nabla v(\cdot,t)\|_{L^q(\Omega)} \le C \quad\mbox{for all}\ t \in (0,T_{\rm max}).
  \]
\end{prop}
%
Unfortunately, the proof of this proposition contains an error in the case $\gamma \in \left(0,\frac{1}{n}\right)$: specifically, in \cite[(3.1) in Section 3]{TZD19} the authors claim that for any $0<\gamma\leq 1$ it is possible to find $s \in \left[1,\frac{n}{(n\gamma-1)_{+}}\right)$ such that 
  \begin{align}\label{s}
    \gamma - \frac{1}{n}<\frac{1}{s}<1 + \frac{1}{n} -\alpha - \beta.
  \end{align}
If from the one hand for $\gamma \in \left[\frac{1}{n},1\right]$ such a relation and \eqref{TZDcondi} fit, from the other hand they do not when $\gamma \in \left(0,\frac{1}{n}\right)$,  and some counterexamples of \eqref{s} can be encountered. For instance, the triplet $(\alpha,\beta,\gamma)=\left(1,\frac{1}{n},\frac{1}{2n}\right)$  is adjusted to \eqref{TZDcondi}, but oppositely it  implies that \eqref{s} is rewritten as $-\frac{1}{2n}<\frac{1}{s}<0$, not satisfied for any $s \ge 1$.
%
%
%
Since relation \eqref{s} is crucial in the derivation of Proposition \ref{prop}, the  machinery to show \cite[Theorem 1.1]{TZD19}, of the item \eqref{thm} above,  misses its validity for $\gamma \in (0,\frac{1}{n})$.
\section{Correction of Proposition \ref{prop} and proof of Theorem \ref{Mainthm}: some hints }
%
%
As specified, we can only confine to the case $\gamma \in \left(0,\frac{1}{n}\right)$. By putting $\gamma_0:=\frac{1}{n}$, we note from \eqref{newcondi} that 
  \begin{align}\label{TZDcondi-gamma0}
    \alpha + \beta + \gamma_0 < 1 + \frac{2}{n}.
  \end{align}
Hence we can fix $s \in [1,\infty)$, rigorously $s \in \left(\frac{1}{\gamma_0},\infty\right)$ (see Remark \ref{remark} below), such that 
  \begin{align}\label{s-gamma0}
    0 = \gamma_0 - \frac{1}{n}<\frac{1}{s}<1 + \frac{1}{n} -\alpha - \beta.
  \end{align} 
We next pick $p \ge \overline{p}$ and $q \ge \overline{q}$, 
where $\overline{p}$ and $\overline{q}$ are defined as in \cite[Section 3]{TZD19}, and set
  \[
    \phi(z):=\int^{z}_{0}\int^{\rho}_{0} \frac{(1+\sigma)^{p-\alpha-2}}{D(\sigma)}\, d\sigma d\rho \quad\mbox{for}\ z \ge 0.
  \]
We can derive (3.9) in \cite{TZD19} unconditionally, that is, we can find $C_1=C_1(q)>0$ such that on $(0,T_{\rm{max}})$ the local solution of problem \eqref{P} complies with
  \begin{align}\label{3.9}
    &\frac{1}{q}\frac{d}{dt} \int_\Omega |\nabla v|^{2q}\, dx 
    + \frac{q-1}{q^2} \int_\Omega |\nabla|\nabla v|^q|^2\, dx
  \\ \notag
  &\le g_1^2\left(2(q-1)+\frac{n}{2}\right) \int_\Omega u^{2\gamma}|\nabla v|^{2(q-1)}\, dx
         + (C_1-2) \int_\Omega |\nabla v|^{2q}\, dx. 
  \end{align}
From the condition $\gamma<\gamma_0$ and Young's inequality it follows that for all $t\in (0,T_{\rm{max}})$ 
  \begin{align*}
    \int_\Omega u^{2\gamma}|\nabla v|^{2(q-1)}\, dx 
    &\le \frac{\gamma}{\gamma_0} \int_\Omega u^{2\gamma_0}|\nabla v|^{2(q-1)}\, dx
         + \left(1-\frac{\gamma}{\gamma_0}\right) \int_\Omega |\nabla v|^{2(q-1)}\, dx
    \\ \notag
    &\le \frac{\gamma}{\gamma_0} \int_\Omega u^{2\gamma_0}|\nabla v|^{2(q-1)}\, dx
           + \left(1-\frac{\gamma}{\gamma_0}\right) \left[\left(1-\frac{1}{q}\right)\int_\Omega |\nabla v|^{2q}\, dx 
           + \frac{|\Omega|}{q}\right].
  \end{align*}
Therefore, by plugging this inequality into \eqref{3.9}, we see that there exist $C_2=C_2(q)>0$ and $C_3=C_3(q,|\Omega|)>0$ providing
  \begin{align}\label{ineq}
    &\frac{1}{q}\frac{d}{dt} \int_\Omega |\nabla v|^{2q}\, dx 
    + \frac{q-1}{q^2} \int_\Omega |\nabla|\nabla v|^q|^2\, dx
  \\ \notag
  &\le C_2\int_\Omega u^{2\gamma_0}|\nabla v|^{2(q-1)}\, dx
         + C_2\int_\Omega |\nabla v|^{2q}\, dx + C_3 \quad \textrm{on} \quad (0,T_{\rm{max}}).
  \end{align}
Since $\gamma_0=\frac{1}{n} \in \left[\frac{1}{n},1\right]$ and \eqref{TZDcondi-gamma0} holds, we can estimate the first term on the right-hand side of \eqref{ineq} as in the proof of \cite{TZD19}, so arriving at \cite[(3.19)]{TZD19}, with  $C_{11}$ involving also the constant $C_3$.
Finally, thanks to relation \eqref{s-gamma0}, we complete the proof by similar arguments to those employed in \cite[Proposition 3.1]{TZD19}.
\qed
%
\begin{remark}[{Comparison between \cite[Theorem 1.1]{TZD19} and \cite[Theorem 2.1]{FV21}}]\label{remark}
The proof of \cite[Proposition 3.1]{TZD19} relies, \textit{inter alia}, on the conservation of mass property $\|u(\cdot,t)\|_{L^{1}(\Omega)}=\int_\Omega u_0(x)\,dx=m$ for all $t\in (0,T_{\rm{max}})$, as well as on the bound  $\|v(\cdot,t)\|_{W^{1,s}(\Omega)} \le C$, valid for any $s \in \left(\frac{1}{\gamma},\frac{n}{(n\gamma-1)_{+}}\right)$, throughout all $t \in (0,T_{\rm max})$ and  for some $C=C(s,\gamma)>0$. The first is obtainable by integrating over $\Omega$ the equation for $u$ in \eqref{P}. For the second, Neumann semigroup estimates, in conjunction with $\int_\Omega g(u)^\frac{1}{\gamma}\leq g_1^\frac{1}{\gamma}m$, entail for some $C_0>0$, $\mu>0$, and all $t\in(0,T_{\rm{max}})$ and $\frac{1}{2}<\rho<1$
\[
\|v(\cdot,t)\|_{W^{1,s}(\Omega)} 
\le C_0\|v_0\|_{W^{1,s}(\Omega)} 
+ C_0 \int^{t}_{0}(t-r)^{-\rho-\frac{n}{2}\left(\gamma-\frac{1}{s}\right)}
e^{-\mu (t-r)}\|u^\gamma(\cdot,r)\|_{L^{\frac{1}{\gamma}}(\Omega)}\, dr.
\]
Conversely, in \cite[Lemma 3.1]{FV21} only a uniform bound for $v(\cdot,t)$ in $W^{1,n}(\Omega)$ and for any $0<\gamma<\frac{2}{n}$ is derived. Subsequently, since $\frac{n}{(n\gamma-1)_{+}}>n$, one concludes that for $s$ close enough to $\frac{n}{(n\gamma-1)_{+}}$, the succeeding $W^{1,s}$-estimates involving $v$, have to play a sharper role on the final result than the $W^{1,n}$-estimates do. This is reflected on condition  \eqref{TZDcondi}, milder than \eqref{FraVigcondi}.
%
\end{remark}
%
 
\end{document}